\input amstex
\magnification=\magstep1
\documentstyle{amsppt}
\def\bC{\bold C}
\def\bR{\bold R}

\def\a{\alpha}
\def\G{\Gamma}
\def\Jac{\hbox{Jac}}

\def\d{\delta}
\def\z{\zeta}
\def\e{\epsilon}
\def\p{\partial}

\def\-{\overline}

\def\O{\Omega}

\def\r{\rho}

\def\r{\rho}
\def\f{\phi}

\def\sm{C^\infty}
\def\cn{\bC^n} 
\def\lt{\mathinner{\raise2pt\hbox{$<$} \mkern-14mu\lower3pt\hbox{$\sim$}}} 

\def\:{\colon}
\def\1p{\O_1^+}
\def\1m{\O_1^-}
\def\2p{\O_2^+}
\def\2m{\O_2^-}
\title 
A CHARACTERIZATION OF THE FINITE MULTIPLICITY\\ OF A CR MAPPING\\ 
\endtitle
\author Yifei Pan\footnotemark\endauthor 
\footnotetext{Department of Mathematics, Indiana-Purdue University Fort Wayne, IN 46805, USA.
 The author was supported in part by a grant
 from Purdue Research Foundation}
\endtopmatter
\heading 1. Introduction\endheading 
In this paper, we give a characterization of the finite multiplicity 
of a CR mapping between real analytic hypersurfaces. 
The finite multiplicity of a CR mapping
was defined algebraically  by Baouendi and Rothschild 
in [BR1] (see the definition below). 
We will
prove that under certain conditions on hypersurfaces the finite 
multiplicity of a CR mapping  is equivalent to that 
the preimage of the map is finite. More precisely,
\proclaim{ Theorem 1} Suppose that $M_1, M_2$ are real analytic 
hypersurfaces of essential finite type in $\cn$ and further $M_2$ contains 
no complex variety 
of positive dimension. Then a smooth CR mapping $f\:M_1\to M_2$ is 
of finite multiplicity at $z_0\in M_1$ is and only if $f^{-1}(f(z_0))$ 
is finite.
\endproclaim
The proof of  Theorem 1 relies on the real analyticity result of[BR1] 
and the
following Theorem 2 that we shall prove. In [BR1], 
Baouendi and Rothschild proved that a smooth CR mapping
of finite multiplicity from a real analytic hypersurface 
of essential finite type
to another real analytic hypersurface is real analytic. 
This result with the proof of Theorem 1 implies 
the following.
\proclaim{Corollary 1} A smooth CR mapping of finite multiplicity 
between real analytic hypersurfaces of essential finite type is the 
restriction of a locally proper holomorphic mapping in $\cn$. 
\endproclaim
\proclaim{Theorem 2} Suppose that $f\:M_1\to M_2$ is a smooth 
CR mapping between real analytic hypersurfaces in $\cn$. Suppose 
further that
$M_1$ is essentially finite and $M_2$ contains no complex variety 
of positive dimension.
If $f^{-1}(f(z_0))\setminus\{z_0\}$ is discrete for a point $z_0\in M_1$,
 then $f$ extends 
holomorphically to a  
neighborhood of $z_0$ in $\cn$.
\endproclaim
A simple example shows that the condition that $M_2$ contains 
no complex variety of positive dimension is necessary in Theorem 1 and 2.
\proclaim{Corollary 2} Suppose that $f\:M_1\to M_2$ is a smooth CR mapping 
between real analytic hypersurfaces of finite type of 
D'Angelo in $\bC^n$. If 
$f^{-1}(f(z_0))\setminus\{z_0\}$ is discrete for a point $z_0\in M_1$
, then $f$ extends holomorphically to a 
neighborhood of $z_0$ in $\bC^n$.
\endproclaim
A well-known problem in the study of real analyticity of CR mappings is 
whether
every smooth CR mapping between real analytic
hypersurfaces of finite type of D'Angelo in $\cn$ is real
analytic. 
\proclaim{Corollary 3} Let $f\:M_1\to M_2$ is a smooth CR mapping 
between real analytic hypersurfaces of finite type of D'Angelo in $\cn$. If
$f$ is real analytic on $M_1\setminus\{p\}$, then $f$ is
also real analytic at $p$.
\endproclaim
This can be reviewed as a "Removable Singularity  Theorem " for the real 
analyticity of CR mappings.
As another corollary of the proof of Theorem 2, one has the following. 
\proclaim{Corollary 4} A 
finite to one smooth CR mapping
from a real analytic hypersurface of essential finite type to another 
real analytic hypersurface is real analytic.
\endproclaim 
Here  a map $f\:M_1\to M_2$ is said to be finite to one if $f^{-1}(q)$ 
is finite for any $q\in M_2$.
The proofs of these results depend  on the work of 
Baouendi-Rothschild [BR1] and Diederich-Fornaess [DF] on real analyticity, 
the Hopf Lemma of [BR3] and  the work of Tumonov [T]
on holomorphic extension of CR functions.
However,
we will directly prove the holomorphic extension of CR mappings whenever 
their work does not apply.  For eariler results, see [L],[Pi],[BJT],[DW], 
[B],[BB] and [BBR]. 
Theorem 1 will be proved in Section 2 and Theorem 2 along with its 
corollaries in Section 3.
The work of this paper is in part inspired by a paper of Pincuk [Pi2]. 
\heading 2.  Proof of Theorem 1\endheading  
To prove Theorem 1, we first recall some basic definitions. Let $M$ be
a real analytic hypersurface in $\cn$ containing the origin and defined 
locally by $\r(z,\-z)=0$, $\nabla\r\not=0$, $z\in\cn$, where $\r$ is a 
real valued analytic function, $\r(0)=0$. As introduced in [BJT], $M$ is
said to be essentially finite at $0$ if for any sufficiently  small
$z\in \cn\setminus\{0\}$ there exists an arbitrarily small $\z\in\cn$ 
satisfying $\r(z,\z)\not=0$, $\r(0,\z)=0$. We point out that if $M$ does 
not contain
any complex variety of positive dimension through $0$, then $M$ is 
essentially
finite at $0$. Consequently, a real analytic hypersurface of 
finite type of D'Angelo is essentially finite.
The finite multiplicity of a CR mapping is introduced by Baouendi and 
Rothschild
in [BR1] as follows. If $f\:M_1\to M_2$ is a smooth CR mapping between 
two smooth real analytic hypersurfaces in $\cn$, there exist $n$ CR 
functions $f_1,...,f_n$ defined on $M_1$
such that $f=(f_1,...,f_n)$. On the other hand if $j$ is a smooth 
CR function defined on $M_1$ near $0$, there exists a formal holomorphic 
power series $J(Z)=\sum 
a_\a Z^\a$ in $n$ indeterminates, such that $U\in u\to Z(u)\in\cn$ ($U$ 
an open
neighborhood of $0$ in $\bR^{2n-1}$, $Z(0)=0$) is a parametrization of $M_1$, 
then 
the Taylor series of $j(Z(u))$ at $0$ is given by $J(Z(u))$. We can choose 
holomorphic coordinates $Z$ such that $\r(Z, 0)=\a(Z)Z_n, \a(0)\not=0$. 
With $Z=(z',z_n)$
and $z'=(z_1,...,z_{n-1})$, the
mapping $f$ is said to be of finite multiplicity at $0$ if
$$\hbox{dim}_{\bC} O[[z']]/(F(z',0))<\infty\tag 1$$
where $(F(z',0)$ is the ideal generated by $F_1(z',0),...,F_n(z',0)$, 
the power series associated to the CR functions $f_1,...,f_n$ and 
$ O[[z']]$ the ring of formal power series in $n-1$ indeterminates 
and the dimension is taken in the sense of vector spaces.

After a holomorphic change of coordinates
 near $0$,  we may assume that $M_1$ is given by an equation
$$\Im z_n=\psi(z',\-z', \Re z_n),\quad \psi(0)
=d\psi(0)=0$$
with $(z',z_n)\in \bC^{n-1}\times \bC$.
We assume that $M_2$ is another real 
analytic hypersurface defined by 
$$\Im z_n=\f(z',\-z', \Re z_n),\quad 
\f(0)=d\f(0)=0$$
with $(z',z_n)\in \bC^{n-1}\times\bC$.
Let $f=(f',f_n)$ be a CR map from $M_1$ 
to $M_2$ with $f(0)=0$. We say that $f_n$ is 
the normal component of $f$ and $z_n$  in the normal 
direction at $0$.
\demo{Proof of Theroem 1} We actually prove that if $f\:M_1\to M_2$ 
is a smooth CR mapping of finite multiplicity between real analytic 
hypersurfaces of essential finite type, then $f^{-1}(0)$ is finite.
By Theorem 1 of [BR1], $f$ is real analytic at $0$. Let 
$F=(F_1(z), ..., F_n(z))$ be the holomorphic extension of $f$ 
to $\cn$ near $0$. If $f^{-1}(0)$ is not finite, then $S=F^{-1}(0)$ 
must be a complex variety of positive dimension. By Theorem 4 of [BR2], 
we have 
$${\p F_n\over\p z_n}(0)\not=0.\tag 2$$
We claim that $S$ lies in $M_1$. Indeed, by $(2)$, 
$$\Im F_n(z)-\f(F'(z),\-F'(z),\Re F_n(z))=0$$
defines a real analytic hypersurface in $\cn$ which clearly conincides 
with $M_1$ near the origin where $F'=(F_1,...,F_{n-1})$. This proves 
the claim.
\enddemo

Now we let $S'$ be any complex curve in $S$ parametrized 
by $$z(\z)=(z_1(\z),...,z_n(\z)$$
passing through $0$. We claim that $z_n(\z)\equiv 0$. Indeed, 
in the chosen coordinates above, by Lemma (3.7) of [BR1], we have
$$F_n(z)=z_nG(z).$$
By (2), we see $G(0)\not= 0$. On $S'$, it follows 
$F_n(z(\z))=z_n(\z)G(z(\z))=0$,
which implies $z_n(\z)=0$. Therefore, $F_1(z',0),...,F_n(z',0)$ 
have common zeros near $0$ and hence  the dimension
$$\hbox{dim}_{\bC} O[[z']]/(F_1(z',0),...,F_n(z',0))$$
is infinite, a contradiction to the finite multiplicity of $f$ at $0$.

As proved above, $S$ lies in $M_1$ and hence $f^{-1}(0)=F^{-1}(0)$. 
This means
that $F$ is a locally proper holomorphic mapping, which gives a proof of 
Corollary 1.

Now we prove that under the conditions in Theorem 1 if $f^{-1}(0)$ is finite
then $f$ is of finite multiplicity. Indeed, by Theorem 2, whose proof 
does not depend on Theorem 1, $f$ is real analytic at $0$. As before, let
$F$ be the holomorphic extension of $f$. 
We  notice $F_n(z)\equiv 0$ since $M_2$ contains
no complex variety of positive dimension and by Theorem 4 of [BR2], $f$ is of finite multiplicity at $0$. This could also proved directly. Indeed, By Theorem
4 of [BR2], (2) holds.
As above, this implies that $F^{-1}(0)$ is finite 
and therefore $F$ is locally proper which implies the finite multiplicity
of $f$.

We close this section by an example. Let $M_1=\{\Im z_3=|z_1|^2+|z_2|^2\}$  
and $M_2=\{\Im z_3=|z_1|^2-|z_2|^2\}$. Consider $f=(g,g,0)$ where is the 
restriction of $e^{-1/z_3^{1/3}}$, which is  holomorphic 
in $\Im z_3>0$ and smooth up the boundary. It is easy to see 
$f^{-1}(0)=0$ but $f$ is not finite multiplicity.  Note that $M_2$ 
contains a complex curve and both $M_1, M_2$
are of essential finite type.
\heading 2.  Proof of Theorem 2\endheading
Following Tumanov [T], we say that a real hypersurface $M_1$ is minimal at
$z_0$ if there is no germ of complex holomorphic hypersurface contained in $M_1$ and passing through $z_0$.
By a theorem of Trepreau [Tr], $f$ extends holomorphically to  one side
of $M_1$. The main result of [BBR] [BR1] [DF] can be stated as 
\proclaim{Theorem} ([BBR][BR1] [DF]) Let $M_1$ is a real analytic 
hypersurface that is essentially finite at $0\in M_1$. If $M_2$ is 
another real analytic hypersurface and $f\:M_1\to M_2$ is a smooth 
CR mapping with $f(0)=0$ and ${\p f_n\over\p z_n}(0)\not=0$, then $f$ 
extends holomorphically to a neighborhood of $0$ in $\cn$.
\endproclaim
The above theorem has
many important applications to global proper holomorphic mappings. 
For example, it was proved in [BR1] [DF] that every proper holomorphic 
mapping between bounded pseudoconvex domains with real analytic boundaries
extends holomorphically across the boundary. 
In [BR2], Baouendi and Rothschild showed that if the normal component
of $f$ is not flat 
(i.e., if there exists a number $k>0$ so that 
${\p^k f_n\over\p z_n^k}(0)\not=0$ )
in the normal direction at $0$ then the condition
$ {\p f_n\over\p z_n}(0)\not=0$ holds automatically. 
As an application of this result, it was proved in [HP]  that the unique 
continuation property holds for proper holomorphic mappings
between bounded domains with real analytic boundaries. This result 
in turn proves that every proper holomorphic mapping between bounded 
real analytic domains that
is smooth up to the boundary extends holomorphically across the boundary.

In order to prove Theorem 2, we need the following lemmas. First we recall 
the definition of a  correspondence. Let $\O$ be a domain in $\cn$
and $f\:\O\to\cn$
be a holomorphic mapping. Denote by $\G_f$ as the graph of 
$f$ $$\G_f=\{(z,w)\:w=f(z), z\in \O\}.$$
Let $$B((z_0,w_0),\e)=\{(z,w)\in\cn\times\cn\:|z-z_0|<\e,|w-w_0|<\e\}.$$ 
We say
that $f$ extends as a correspondence to a neighborhood of $(z_0,w_0)$ if 
there exist $\e>0$ and a pure n-dimensional subvariety 
$$V\subset B((z_0,w_0),\e)$$
such that 
$$\G_f\cap B((z_0,w_0),\e)\subset 
V\cap B((z_0,w_0),\e).$$
Now we state a lemma due to Bedford and Bell [BB].
\proclaim {Lemma 1} Let $\O$ be a bounded domain in $\cn$ with 
smooth boundary
near $z_0\in\p\O$, and let $f\:\O\to\cn$ be a holomorphic mapping that 
is $\sm$ smooth up to the boundary of $\O$ near $z_0$. Then $f$ extends 
holomorphically to a neighborhood of $z_0$ in $\cn$.
\endproclaim
Let $M_1$ and $M_2$ be smooth real hypersurfaces in $\cn$ and let
$\O_1, \O_2$ be two domains in $\cn$ with defining functions $r_i$ 
for $i=1,2$ such that
$\nabla r_i\not=0$ on $\O_i$ for $i=1,2.$
Set $\O_i^+=\{z\in \O_i\:r_i(z)>0\}$ and $\O_i^-=\{z\in \O_i\:r_i(z)<0\}$ 
for $i=1,2.$

If $F\:\1m\to\cn$ is a holomorphic mapping, we denote by $\Jac F$  
the determinant of the Jacobian matrix of $F$.

As will become clear, in order to prove Theorem 2, 
one has to only consider the case when
$w_0$ is a minimal but not minimally convex point in the sense of [BR3]. 
For this matter, we prove the following result.
\proclaim {Lemma 2} Let $f\:M_1\to M_2$ be a smooth CR mapping between 
smooth real hypersurfaces $M_1, M_2$ in $\cn$. Suppose that $f$ extends 
holomorphically to an one-sided neighborhood of $M_1$, say $\1m$. 
Given a point $z_0\in M_1$, if $M_2$ contains no nontrivial complex 
variety through $f(z_0)$ and if $f(z_0)$ is not 
minimally convex and  $f^{-1}(f(z_0))\setminus\{z_0\}$ is discrete   
, then $f$ extends holomorphically to a neighborhood of
$z_0$ in $\cn$.
\endproclaim
We remark that no real analyticity on hypersurfaces is assumed above. 
\demo{Proof of  Lemma 2}  Let $F(z)\:\O_1^-\to\cn$ be the extension of $f$. 
First we prove two facts to be used later.

We notice that $F(\1m)\not\subset M_2$ since $M_2$ contains no complex 
variety of positive dimension.
Now we claim that 
$\Jac F(z)\not\equiv 0$. Indeed, 
if $\Jac F(z)\equiv 0$ in $\1m$, we let $\mu$
be the maximal rank of the Jacobian matrix of $F$ in $\1m$. 
We have $0<\mu<n$ and the set 
$$\{z\in\1m\: \hbox{Rank} F=\mu\}$$ 
is
an open dense subset of $\1m$.
By the rank 
theorem
and the fact $F(\1m)\not\subset M_2$, we may find a sequence of 
points $z_k\in\1m$ converging to $z_0$ such that $F(z_k)\not\in M_2$ 
and for each $k$ the analytic set 
$$\{z\in \1m\:F(z)=F(z_k)\}$$ 
has an irreducible component $V_k\subset \1m$ of
dimension $n-\mu>0$ passing through $z_k$. Since $F(z_k)\not\in M_2$, it 
follows that for each $k$,  $V_k$ does not have limit points on $M_1$. 
Therefore $\-V_k$ is
a closed analytic variety in $\O_1$. Now let 
$z'\in\-{\cup\-V_k}\setminus \cup\-V_k$, and we see that 
$$f(z')=F(z')=\lim F(z_k)=F(z_0)=w_0$$ 
This implies that $z'\in f^{-1}(w_0)$.
But $f^{-1}
(w_0)\setminus\{z_0\}$ is discrete, we  see  that the sequence of 
the sets $\-V_k$ clusters on $M_1$ only at  discrete points near $z_0$. 
Thus by the generalized continuity
principle we conclude that $F(z)$ extends holomorphically to a 
neighborhood of $z_0$ in $\cn$. As before, it implies 
that $\Jac F(z)\not\equiv 0$ since $F$ is locally proper.
\enddemo
Using these facts we will prove that $F$ extends holomorphically to a 
neighborhood of $z_0$ in $\cn$. 

 When $w_0\in M_2$ is  not a minimally
convex point, an important fact  is  that (see Theorem 7 
of [BR4] , Theorem 1 of [BR3], [T])
 every holomorphic function defined
on one side of $M_2$, which admits a distribution
limit up to $M_2$, extends holomorphically to a small open neighborhood
of $w_0$.  This fact has been used in [HP], [P].

To be able to prove the holomorphic extension of $F$ when $w_0$ is 
not minimally convex, we will construct pieces of 
proper holomorphic mappings near $z_0$.

Since $f^{-1}(w_0)\setminus\{z_0\}$ is discrete and $f^{-1}(w_0)$ is closed,  
we may choose an open neighborhood $\O_1$ of $z_0$ such that 
$$\p\O_1\cap \{f^{-1}(w_0)\}=\emptyset.$$
So we have that $\hbox{dist}(\p\O_1, \{f^{-1}(w_0)\})=\d>0$. 

Now consider 
$$V=\{z\in\1m, F(z)=w_0\}.$$
Then $V$ is an analytic variety in $\1m$. If $\hbox{dim} V\geq 1$, 
let $ V'$ be 
an irreducible component of $V$.  Since $V$ only has limit 
points $f^{-1}(w_0)$ on $M_1$, by Shiffman's theorem, $\-V'$ is 
an analytic variety in $\O_1$.
The continuity principle implies that $F$ extends holomoprhically to 
a neighborhood of $z_0$.

Now we may assume that $\hbox{dim} V=0$. This means $V$ is a 
discrete set in $\1m$. We may shrink $\O_1$ slightly so 
that $\p\O_1\cap V=\emptyset$. Therefore
, we have 
$$\hbox{dist}(w_0, F(\p\1m\setminus M_1))>0.$$
Then we can choose a very small open neighborhood $\O_2$ of $w_0$ 
such that $$\hbox{dist}(\p\O_2, F(\p\1m\setminus M_1))>0.\tag \#$$   
Since $F(\1m)\not\subset M_2$, $F(\1m)$ intersects at least one side 
of $M_2$. Therefore there are two possibilities as follows.
(I)   For any small neighborhood $\O_2$ of $w_0$ we 
have 
$$F(\1m)\cap{\O_2^-}\not=\emptyset\quad\hbox{and}
\quad F(\1m)\cap\O_2^{+}\not=\emptyset.$$

(II) There is an arbitrarily small neighborhood $\O_2$ of $w_0$ such that
$$ F(\1m)\subset \-{\O_2^-}\quad\hbox{or}\quad F(\1m)\subset \-{\O_2^+}.$$

We consider the case (I) first, the case (II) can be delt with similarly.

Consider two nonempty open sets in $\1m$: 
$$U^+=F^{-1}(\O_2^+)\quad\hbox{and}\quad U^-=F^{-1}(\O_2^-).$$
We claim that the restriction of $F$ to $U^+$ (resp. $U^-$) is a proper 
map from $U^+$ to $\O_2^+$(resp. $\O_2^-$).
Indeed, let $F^+=F|U^+$ and let $K\subset\subset \O_2^+$ be a compact subset, 
we want to prove that $(F^+)^{-1}(K)$ is a compact subset in  
$U^+$. If $(F^+)^{-1}(K)$ is not compact in $U^+$, there exists a 
point $p\in\p U^+$ such that $F^+(p)\in K$. Since $K\cap M_2=\emptyset$, 
we have $p\not\in M_1$, and by ($\#$) $p\in\1m$. Therefore, 
there exists a neighborhood
$O$ of $p$, such that $F(O)\subset\O_2^+$. Hence $p$ cannot be a boundary 
point of $U^+$, a contradiction.

Now we observe that the open set $U^+\cup U^-$ is, in general, 
not connected. We make some simple observations that are crucial to 
what follows in the proof of Lemma 2.

Claim 1:  The set $U^+\cup U^-$ is an open dense set near $z_0$ in $\1m$ 
along $M_1$.
Indeed, if it is not the case, then, for any small neighborhood of $z_0$, 
there
exists a point $p\in \1m$ in that neighborhood, 
and there exists a small neighborhood
$O$ of $p$ contained in $\1m$ so that $f(O)\subset \p\O_2\cup M_2$
(since by continuity $F(O)\subset\-\O_2$). This
is impossible since $\Jac F(z)\not\equiv 0$ in $\1m$.

Claim 2:  The open set $U^+\cup U^-$ has finitely many connected components.

Indeed, if it is not the case, we let $U_j$ be connected components of 
$U^+\cup U^-$ for $j=1,2,...$. Let $E_j=\p U_j$ be the boundary of $U_j$. 
Since $\1m$ is bounded,  either $\{E_j\}$ 
accumulates at a neighborhood of an interior point of $\1m$ 
where they are disjoint each other, or at a boundary point  or both.
We prove that neither is possible. Indeed, if $E_j$ accumulates at $p\in\1m$ 
we can assume that $\Jac F(p)\not=0$ since $\Jac F\not\equiv 0$ and 
$\{z\in\1m, \Jac F(z)=0\}$ is an analytic variety of complex dimension of 
$n-1$. Therefore $F$ is a  local biholomorphism in a neighborhood $O$ of 
$p$, 
 therefore we may assume $E_j\subset O$ locally near $p$ for all $j$.
On the other hand, we have  $F(E_j)\subset
\p\O_2\cup M_2$ for all $j$, from which we arrive at  a contradiction.

If $E_j$ accumulates at $p\in M_1$ 
we can assume that $\Jac F(p)\not=0$ since $\{z\in M_1\: \Jac F\not=0\}$ 
is a dense open subset of $M_1$.
Then the above argument applies since $F$ is a diffeomorphism near $p$ 
after we extend $F$ smoothly to a neighborhood of $p$ in $\cn$.

Now let $\{U_j^+\}_{j=1}^k$ be connected components of $U^+$, 
similarly $\{U_j^-\}_{j=1}^l$ for $U^-$. 
Let $g_j$ be the restriction of $F$ on $U_j^+$,
and $h_j$ on $U_j^-$. Therefore $g_j\: U_j^+\to\O_2^+$ and
$h_j: U_j^-\to \O_2^-$ are proper holomorphic mappings.

We then consider a proper mapping $g$ from $D$ to $G$, where 
the paring $(D, G)$ is  either $(U_j^+,\O_2^+)$ or  $(U_j^-,\O_2^-)$ 
and $g$ is either $g_j$ or $h_j$.
The graph of $g$ is defined to be
$$\G_g=\{(z,w)\in D\times G, w=g(z)\}.$$ By the Proper Mapping Theorem, 
$g$ is a covering from $D\setminus g^{-1}(g(V_g))$ to $G\setminus g(V_g)$ 
of multiplicity $m$, where 
$$V_g=\{z\in D\: \Jac g=0\}.$$  

Let $G_1, G_2,...,
G_m$ be the local inverses defined on $ G\setminus g(V_g)$. Define 
over $D\times G\setminus g(V_g)$ $$H_i(z,w)=\Pi_{j=1}^m(z_i-(G_j(w))_i).$$
By the removable singularity result of bounded holomorphic functions,  
$H_i$ extends to be holomorphic on $D\times G$.
Denote
$$A_g=\{(z,w)\in D\times G: H_1=H_2=...=H_n=0\}.$$
It is easy to check that $\G_g=A_g$.

Let $\G_{g_j}$, $\G_{h_j}$ be the graphs of $g_j, h_j$ respectively,
and let $A_{g_j}$, $A_{h_j}$ be associated with $g_j, h_j$ as defined above.
We see that the graph of $F$ over $U^+\cup U^-$ is given
by 
$$\cup_{j=1}^k\G_{g_j}\cup \cup_{j=1}^l\G_{h_j},$$ 
which is equal to
$$\cup_{j=1}^k A_{g_j}\cup \cup_{j=1}^l A_{h_j}.$$  
As we have observed that
the open set $U^+\cup U^-$ is an open dense set along $M_1$ near $z_0$ 
(Claim 1). By the continuity, we conclude that the graph of $F$ over a 
small one-sided neighborhood of $M_1$ near $z_0$ is contained in 
$$\cup_{j=1}^k A_{g_j}\cup \cup_{j=1}^l A_{h_j}.$$ 

 Now we want to show that 
 $$\cup_{j=1}^k A_{g_j}\cup \cup_{j=1}^l A_{h_j}$$ 
 extends to be an analytic variety of pure dimension n in $\cn\times\cn$
 near $(z_0,w_0)$.  Indeed, we notice for each $g$ ( either $g_i$ or $h_j$)
 $$H_i(z,w)=z_i^m+S_{m-1}(w) z_i^{m-1}+...+S_0(w).$$
where $S_j(w)$ is  the j-th symmetric function of $(G_j(w))i$ for $j=1,..,m$.
 Since $S_j(w)$ are bounded, and  since $w_0$ is not minimally 
 convex, then $S_j(w)$ extends to be holomorphic  in a neighborhood of
 $w_0$ in $\cn$ from
 either side whenever applicable. 
 Therefore $H_i(z,w)$ extends to be holomorphic to a 
neighborhood of $(z_0,w_0)$ in $\cn\times\cn$, this , in turn, implies that
 $$\cup_{j=1}^k A_{g_j}\cup \cup_{j=1}^l A_{h_j}$$ 
 is an analytic variety of
pure dimension $n$ in a neighborhood of $(z_0, w_0)$, which implies that
 $F$ extends to be a correspondence to  a neighborhood of $z_0$.
 Lemma 1 then gives the holomorphic extension of $F$ at $z_0$. This completes
the proof of Lemma 2 for the case (I). Case (II) can be proved equally.
\demo{Proof of Theorem 2} Let $z_0\in M_1$, $w_0=f(z_0)\in M_2$. 
Since $M_1$ is minimal at $z_0$, by Trepreau's Theorem, $f$ extends 
holomorphically to one side neighborhood of $M_1$, say $\1m$, 
the extension is denoted by $F(z)$. Therefore $F(z)\:\1m\to\cn$ is 
a holomorphic mapping, such that $F=f$ on $M_1$. 
If $w_0$ is minimally convex, then the complex Hopf Lemma
of [BR3] and the theorem of [BR1] and [DF] imply the $f$ extends 
holomorphically to a neighborhood of $z_0$ since $\hbox{Jac} F\not\equiv 0$. 
When $w_0$ is not minimally convex then Lemma 2 applies. This completes 
the proof.
\enddemo

Corollary 2 is a special case of Theorem 2 since a real anlytic 
hypersurface of finite type of D'Angelo is essentially finite and 
contains no nontrivial complex varieties.

Now we give a proof of Corollary 3.
\demo{Proof of Corollary 3}  It suffices  to prove that 
$f^{-1}(f(p))\setminus\{p\}$
is discrete. Let $q \in f^{-1}(f(p))\setminus\{p\}$ but $q\not=p$. 
We want to prove that $q$ is an isolated point.
Since $f$ is real analytic at $q$ by the assumption, then $f$ extends 
holomorphically to a neighborhood of $q$, say the extension as $F$.
By  a result of [BR2], the Hopf Lemma holds at $q$ for the normal 
component of $F$.  Let $\r$ be a real analytic defining function of 
$M_2$ near $w_0$. By the Hopf Lemma just mentioned at $q$, it is easy to see,
by changes of coordinates at both $q$ and $f(q)$, 
that $\r\circ F$ is again a defining function of $M_1$ near $q$. 
Therefore the equation
$$\{z\in\cn\:\r\circ F(z)=0\}$$
defines a real analytic hypersurface near $q$, which is 
identical to $M_1$ near $q$.
This implies that $F^{-1}(f(q))$ is contained in $M_1$. Since $M_1$
is of finite type of D'Angelo and $F^{-1}(f(q))$ is 
a complex analytic variety, we conclude that $q$ is a isolated
point in $M_1$. Theorem 2 then applies at $p$ since 
$f^{-1}(f(p))\setminus \{p\}$
is discrete.

In order to prove Corollary 4, we prove the following first.
\proclaim{Lemma 3} Let $f\: M_1\to M_2$ be a  finite to one  
smooth CR mapping
between smooth real hypersurfaces that extends holomorphically to
$\1m $ as $F$.
Given $z_0\in M_1$ and $f(z_0)=w_0\in M_2$.
If $M_2$  is minimal but not minimally convex at $w_0$,
then $f$ extends holomorphically  to a 
neighborhood of $z_0$ in $\cn$
\endproclaim
\demo{Proof}  By the proof of Lemma 2, it suffices to prove the following
 (I) $F(\1m)\not\subset M_2$, and  
(II) $\Jac f(z)\not\equiv 0.$ 

Indeed,  if 
$F(\1m)\subset M_2$ then $\Jac F(z)\equiv 0$ in $\1m$. This implies 
that the Jacobian matrix of the map $f:M_1\to M_2$ considered as a 
real map of the real manifolds is of maximal rank $\mu$ such that 
$0<\mu<2n-1$. Therefore by the rank theorem, there exists a point $w'$ 
near $w_0$ 
such that $f^{-1}(w')$ is a manifold of dimension $n-\mu$, a 
contradiction to finite to one.  (II) follows too.
\enddemo
\demo{Proof of Corollary 4} First we observe that $F(\1m)\not\subset M_2$ 
by the proof of Lemma 3. If $w_0$ is not minimal, then, 
by a unique continuation result for holomorphic mappings in (Theorem 2, [P]), 
$F$ does not vanish
to infinite order at $z_0$ in the normal component. Then $F$ extends 
holomorphically to a neighborhood of $z_0$. The rest of the proof 
follows as in Theorem 2 by Lemma 3.
\enddemo
\end